\documentclass[reqno,final]{amsart}
\usepackage{natbib}  
\usepackage{fancyhdr} 
\usepackage{color} 
\usepackage{hyperref} 
\usepackage{graphicx} 

\definecolor{aleacolor}{rgb}{0.16,0.59,0.78}

\hypersetup{
breaklinks,
colorlinks=true,
linkcolor=aleacolor,
urlcolor=aleacolor,
citecolor=aleacolor}

\renewcommand{\headheight}{11pt}
\renewcommand{\cite}{\citet}

\theoremstyle{plain}
\newtheorem{theorem}{Theorem}[section]

\theoremstyle{definition}

\theoremstyle{remark}
\newtheorem{remark}[theorem]{Remark}

\makeatletter \@addtoreset{equation}{section} \makeatother
\renewcommand\theequation{\thesection.\arabic{equation}}
\renewcommand\thefigure{\thesection.\arabic{figure}}
\renewcommand\thetable{\thesection.\arabic{table}}

\newcommand{\aleaIndex}[1]{\href{http://alea.impa.br/english/index_v#1.htm}{\bf #1}}
\eheader{Alea}{\aleaIndex{7}}{2010}{193}{205}

\elogo{\framebox[5.8cm][c]{\footnotesize \parbox[c]{5.5cm}{
The original article is published by the \href{http://alea.impa.br/english/index_v7.htm}{Latin American Journal of Probability and Mathematical Statistics}
} } \parbox[c]{3cm}{\includegraphics[width=3cm]{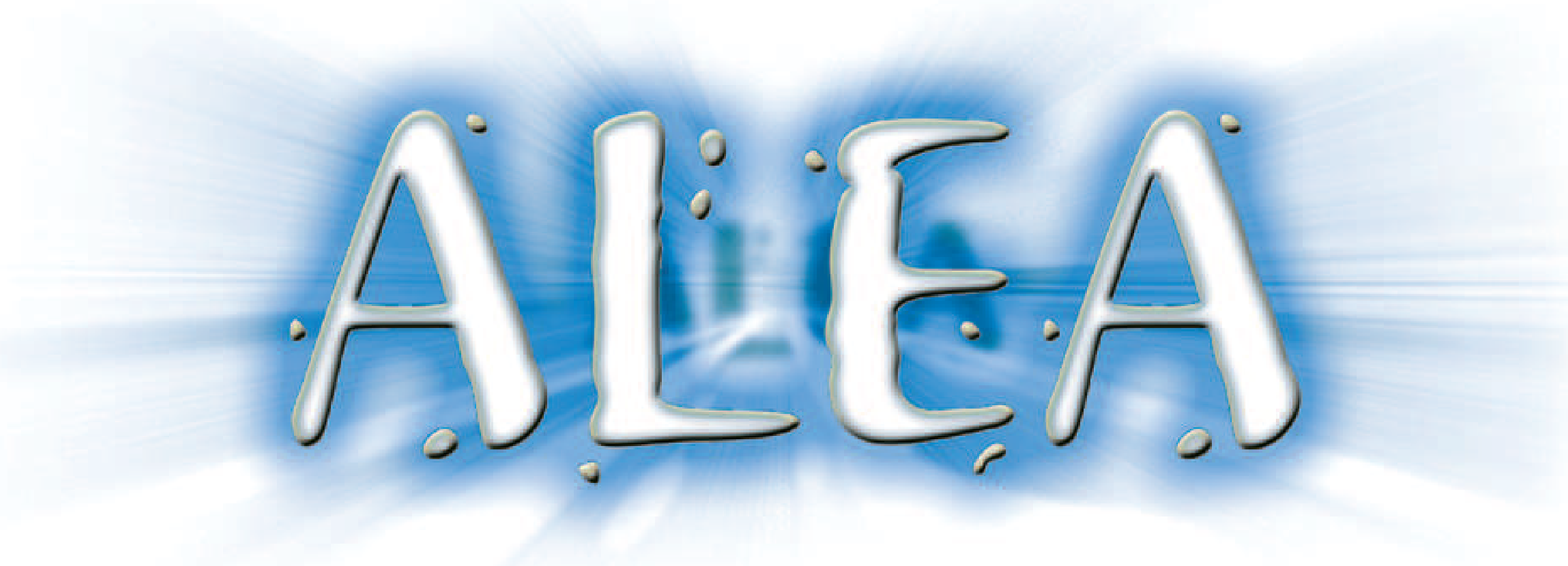}}}

\usepackage{amsmath}
\usepackage{amsthm}
\usepackage{amssymb}
\usepackage{latexsym}
\usepackage{mhequ}

\def\eps{\varepsilon}

\def\Z{\mathbb {Z}}

\def\A{{\mathcal A}}

\def\bd{{\bf d}}
\def\wh{\widehat}
\def\sign{\operatorname{sign}}

\def\vphi{\varphi}
\def\bone{{\bf 1}}

\def\ep{\varepsilon}

\def\si{\sigma}

\newcommand{\comment}[1]{}

\newcommand{\eq}{\begin{equation}}
\newcommand{\en}{\end{equation}}

\def\qed{{\hfill $\Box$ \bigskip}}

\begin{document}

\title{Crowding of Brownian spheres}
\date{February 4, 2010; accepted June 1st, 2010}

\author{Krzysztof Burdzy}
\address{Department of Mathematics, Box 354350,
University of Washington, Seattle, WA 98195}
\address{http://www.math.washington.edu/\~{}burdzy/ }
\email{burdzy@math.washington.edu}
\author{Soumik Pal}
\address{Department of Mathematics, Box 354350,
University of Washington, Seattle, WA 98195}
\address{http://www.math.washington.edu/\~{}soumik/}
\email{soumik@math.washington.edu}
\author{Jason Swanson}
\address{Department of Mathematics, University of
Central Florida, 4000 Central Florida Blvd.,
P.O. Box 161364, Orlando, FL 32816-1364}
\address{http://www.swansonsite.com/W/}
\email{swanson@mail.ucf.edu}

\thanks{K.B.: Research supported in part by NSF Grant DMS-0906743 and
by grant N N201 397137, MNiSW, Poland.}
\thanks{J.S.: Research supported in part by NSA Grant H98230-09-1-0079.}

\begin{abstract}
We study two models consisting of reflecting one-dimensional
Brownian ``particles'' of positive radius. We show that the
stationary empirical distributions for the particle systems do
not converge to the harmonic function for the generator of the
individual particle process, unlike in the case when the
particles are infinitely small.
\end{abstract}


\maketitle

\section{Introduction }

In this article we consider the dynamics of a collection of
\textit{hard Brownian spheres} with drifts or boundary
conditions that includes instantaneous reflections upon
collisions. The models are similar to existing ones in the
literature that consider point masses instead of spheres of a
positive radius. We will show that the (empirical) distribution
of a family of Brownian spheres behaves differently from the
(empirical) distribution of the point Brownian particles in
some natural models. In particular, the distribution of
Brownian spheres fails to satisfy the usual heat equation under
circumstances that lead to the heat equation for the infinitely
many infinitesimally small Brownian particles.

Various models of colliding Brownian particles have been
considered in the statistical physics literature. One stream,
pioneered by \cite{H}, considers a countable collection
of Brownian point masses on the line that collide and reflect
instantaneously. Also see the follow-up work on tagged particle
in the Harris model by D\"urr, Goldstein, \& Lebowitz,
\cite{DGL}. A variation on the theme has been to replace the
instantaneous reflection by a potential, and goes by the name
of \textit{gradient systems}. In these gradient systems,  one
studies the behavior of countably many particles under a
repelling potential. Usually the potential is modeled as smooth
with a singularity at zero; see the article by \cite{CL}. 
A particular example of this class
includes the famous Dyson Brownian motion from Random Matrix
theory; see \cite{D62}, and 
\cite{CL2}. The other class of models, closer to our article,
goes by the name of \textit{hard-core} interactions, in which
the Brownian particles are assumed to be hard balls of small
radius, and consequently, there is instantaneous reflection
whenever two such balls collide (plus possible additional
interactions). This is the spirit taken in the articles by
Dobrushin \& Fritz in dimension one, \cite{DF}, and Fritz \&
Dobrushin in dimension two, \cite{FD}, \cite{L77,L78} (with a correction by \citealp{S79}). The main
focus of these authors is the non-equilibrium dynamics of the
gradient systems. Also see the articles by \cite{O96,O98}, and \cite{T96} all of which consider
properties of a tagged particle in the infinite system.

In the discrete case, the various models of symmetric and
asymmetric exclusion processes have been considered. Closest in
spirit to the models discussed here is the totally asymmetric
exclusion process (TASEP) considered by Baik, Deift, \&
Johansson, \cite{BDJ}, and \cite{J00} in connection
with random matrices and the longest increasing subsequence
problem. Specifically, if the initial configuration in TASEP is
$\mathbb{Z}^-$, then the probability that a particle initially
at $-m$ moves at least $n$ steps to the right in time $t$
equals the probability distribution of the largest eigenvalue
in a unitary Laguerre random matrix ensemble. In recent
subsequent articles \cite{TW,TW09}, Tracy \& Widom
explicitly compute transition probabilities of individual
particles in the asymmetric exclusion process, extending
Johansson's work.

In this paper, we consider only one dimensional models, so our
``spheres'' are actually intervals. The title of this paper
reflects our intention to study multidimensional models in
future articles.  We leave more detailed discussion to Section~\ref{sec:disc}. That section also contains references to
related research projects.

We consider two models, which have the following common
features. Informally speaking, both models consist of families
of Brownian ``particles''. The $k$-th ``particle'' is
represented by an interval $I^k_t = (X^k_t-\eps/2,
X^k_t+\eps/2)$, where $X^k_t$ is a Brownian-like process. The
intervals $I^k$ and $I^j$ are always disjoint, for $k\ne j$.
The processes $X^k$ are driven by independent Brownian motions.
When two intervals $I^k$ and $I^j$ collide, they reflect
instantaneously. In the first model, the number of particles is
constant and they are pushed by a barrier moving at a constant
speed. In the second model, particles enter the interval
$[0,1]$ at the left, they reflect at 0, and they are killed
when they hit the right endpoint. The second model is our
primary focus because it is related to other models considered
in mathematical physics literature---see Section
\ref{sec:disc}.

We are grateful to Thierry Bodineau, Pablo Ferrari, Claudio
Landim, Mario Primicerio and Jeremy Quastel for very helpful
advice. We would like to thank the referee for the suggestions
for improvement.

\section{Extreme crowding}

We start with an informal description of our first model, which
consists of a fixed number $n$ of ``particles''. The $k$-th
leftmost ``particle'' is represented by an interval $I^k_t =
(X^k_t-\eps/2, X^k_t+\eps/2)$. The intervals $I^k$ and $I^j$
are always disjoint. The processes $X^k$ are driven by
independent Brownian motions. When two intervals $I^k$ and
$I^j$ collide, they reflect instantaneously. The intervals are
pushed from the left by a barrier with a constant velocity,
that is, the leftmost interval reflects on the line $x= ct$.

Formally, we define $\{X^0,X^1,\ldots,X^n\}$ to be continuous
processes such that $X^0_t=-\ep/2+ct$, $X^k_t-X^{k-1}_t\ge\ep$
for all $k\ge1$ and all $t\ge0$, and for $k\ge 1$,
  \[
  dX^k_t = dB^k_t + dL^k_t - dM^k_t,
  \]
where $\{B^1,\ldots,B^n\}$ are iid Brownian motions, and $L^k$
and $M^k$ are nondecreasing processes such that
  \[
  \int_0^\infty 1_{\{X^k_t-X^{k-1}_t>\ep\}}\,dL^k_t = 0
    \qquad\text{and}\qquad
    \int_0^\infty 1_{\{X^{k+1}_t-X^k_t>\ep\}}\,dM^k_t = 0.
  \]
(Here, we may interpret $X^{n+1}\equiv\infty$.) The
distributions of $X^k_0$ for $1\le k\le n$ will be specified
later.

To construct the solution to this Skorohod problem, consider
first the processes $Y^k_t=X^k_t-(k-1)\ep-\ep/2-ct$. These
processes satisfy $Y^0\equiv0$, $Y^k_t-Y^{k-1}_t\ge0$ for all
$k\ge1$ and all $t\ge0$, and for $k\ge 1$, $dY^k_t = dB^k_t -
c\,dt + dL^k_t - dM^k_t$, where
  \[
  \int_0^\infty 1_{\{Y^k_t-Y^{k-1}_t>0\}}\,dL^k_t = 0
    \qquad\text{and}\qquad
    \int_0^\infty 1_{\{Y^{k+1}_t-Y^k_t>0\}}\,dM^k_t = 0.
  \]
We may therefore construct the processes $\{Y^1,\ldots,Y^n\}$
using order statistics.\break Namely, let $\{Z^1,\ldots,Z^n\}$ be
defined by $dZ^k_t = dB^k_t - c dt$, and reflected at 0. For
every fixed $t\geq 0$, we let $Y^1_t, Y^2_t, \dots, Y^n_t$ be
ordered $Z^k_t$'s, that is, $\{Y^1_t, \dots, Y^n_t\} =\{Z^1_t,
\dots, Z^n_t\}$ and $Y^1_t \leq Y^2_t\leq \dots \leq Y^n_t$.
Finally, we let $X^k_t = Y^k_t +(k-1) \eps + \eps/2 + c t$ and
$I^k_t = (X^k_t-\eps/2, X^k_t+\eps/2)$.

Let $n\eps = b$. We will fix $b>0$ and analyze the behavior of
the system of intervals $\{I^k\}$ as $n\to\infty$. In other
words, we will keep the total length of all intervals $I^k$
constant.

The stationary distribution for $Z^k$ has the density $\vphi(z)
= ce^{-cz}$ for $z\geq 0$, with $c= 2c_1$, because
\begin{align*}
\frac 1 2 \frac{d^2}{dz^2} \vphi(z) + c_1 \frac d {dz} \vphi(z) =0.
\end{align*}

Consider any $0\leq x_1 < x_2 < \infty$, let $\lambda$ denote the
Lebesgue measure, and let
  \begin{align}\label{eq:m17.4}
  \bd([x_1, x_2]) = \bd_t([x_1, x_2])
    = \frac { \lambda\left( [x_1 + ct, x_2 + ct]
    \cap \bigcup_{1\leq k \leq n} I^k_t  \right)}
    {x_2 - x_1}.
  \end{align}
The quantity $\bd([x_1, x_2])$ represents the average density
of ``particles'' $I^k$ on the interval $[x_1, x_2]$.

We will say that the intervals $\{I^k\}$ have the
pseudo-stationary distribution if all $Z^k_t$'s are independent
and have the stationary distribution $\vphi$ for $t=0$ and,
therefore, for every $t\geq 0$.

\begin{theorem}\label{th:m16.1}
Suppose that the intervals $\{I^k\}$ have the pseudo-stationary
distribution. Fix arbitrary $p_1, d_1< 1$, $d_2>0$, and $0\leq
x_1 < x_2 < b < x_3 < x_4 < \infty$. There exist $c_0, n_0 <
\infty$ such that for $c\geq c_0$, $n \geq n_0$ and any $t\geq
0$, we have
\begin{align}\label{eq:m17.1}
P(\bd_t([x_1, x_2]) \geq d_1) &\geq p_1,\\
P(\bd_t([x_3, x_4]) \leq d_2) &\geq p_1. \label{eq:m17.2}
\end{align}
\end{theorem}

The theorem says that the ``particles'' $I^k$ clump together
and there is a sharp transition in density of ``mass'' around
$x = b$. This is in contrast with infinitely small
``particles'' $Z^k$ whose empirical distribution is close to
the distribution with the density $\vphi(z) = c e^{-cz}$ that
displays no sharp drop-off.

\begin{proof}[Proof of Theorem \ref{th:m16.1}]
Without loss of generality, we let $t=0$. We define
$y_k\in(0,\infty)$ in an implicit way by the following formula,
for $k=1, 2,3,4$,
\begin{align*}
x_k = b \int_0^{y_k} \vphi(z) dz + y_k.
\end{align*}

Note that $y_1 < y_2$, and that for $\ep>0$ sufficiently small
(that is, for $n=b\ep^{-1}$ sufficiently large), it is possible
to choose $y_5,y_6$ such that $y_1< y_5 < y_6 < y_2$, and
\begin{align}\label{eq:m16.2}
\frac
{b \int_{y_5}^{y_6} \vphi(z) dz- 2\eps}
{y_2 - y_1 + b \int_{y_1}^{y_2} \vphi(z) dz  }
\geq
\frac
{b \int_{y_1}^{y_2} \vphi(z) dz}
{y_2 - y_1 + b \int_{y_1}^{y_2} \vphi(z) dz  }
- (1-d_1)/2.
\end{align}
Since $b-x_2 >0$, we can find $c$ so large that,
\begin{align}\label{eq:m16.4}
\frac
{c (b -x_2 ) }
{1 + c (b -x_2 )  }
\geq 1- (1-d_1)/2.
\end{align}

Let $\lceil a\rceil$ denote the smallest integer greater than
or equal to $a$. By the law of large numbers, if $n$ is
sufficiently large, the number of $Z^k_0$'s in the interval
$[0,y_1]$ is smaller than or equal to $n\int_{0}^{y_5} \vphi(z)
dz$, with probability greater than $1-(1-p_1)/2$. If this event
holds then there are exactly $\left\lceil n\int_{0}^{y_5}
\vphi(z) dz\right\rceil$ processes $Z^k_0$ in some (random)
interval $[0,y_7]$ with $y_7 \geq y_1$. This implies that there
are exactly $\left\lceil n\int_{0}^{y_5} \vphi(z)
dz\right\rceil$ processes $X^k_0$ in $[0, \eps \left\lceil
n\int_{0}^{y_5} \vphi(z) dz\right\rceil + y_7]$. Note that
\begin{align*}
\eps \left\lceil n\int_{0}^{y_5} \vphi(z) dz\right\rceil + y_7
\geq b \int_{0}^{y_5} \vphi(z) dz + y_7
\geq b\int_{0}^{y_1} \vphi(z) dz + y_1 = x_1.
\end{align*}
Hence, the number of $X^k_0$'s in the interval $[0, x_1]$ is
smaller than or equal to $n\int_{0}^{y_5} \vphi(z) dz+1$, with
probability greater than $1-(1-p_1)/2$. A completely analogous
argument shows that, if $n$ is sufficiently large, then the
number of $X^k_0$'s in the interval $[x_2, \infty]$ is smaller
than or equal to $n\int_{y_6}^{\infty} \vphi(z) dz+1$, with
probability greater than $1-(1-p_1)/2$. Both events hold with
probability greater than $1-2(1-p_1)/2 = p_1$, and then the
number of $X^k_0$'s in $[x_1, x_2]$ is greater than or equal to
$n\int_{y_5}^{y_6} \vphi(z) dz-2$. This and \eqref{eq:m16.2}
imply that
\begin{align}\label{eq:m16.3}
\bd([x_1,x_2]) &\geq \frac
{\eps n\int_{y_5}^{y_6} \vphi(z) dz-2\eps}
{x_2 - x_1 }
= \frac
{b \int_{y_5}^{y_6} \vphi(z) dz - 2\eps}
{b \int_0^{y_2} \vphi(z) dz + y_2 - b \int_0^{y_1} \vphi(z) dz - y_1 }
\nonumber \\
&= \frac
{b \int_{y_5}^{y_6} \vphi(z) dz - 2\eps}
{y_2 - y_1 + b \int_{y_1}^{y_2} \vphi(z) dz  }
\geq
\frac
{b \int_{y_1}^{y_2} \vphi(z) dz}
{y_2 - y_1 + b \int_{y_1}^{y_2} \vphi(z) dz  }
- (1-d_1)/2.
\end{align}

We have
\begin{align*}
x_2 = b \int_0^{y_2} \vphi(z) dz + y_2
= b \int_0^{y_2} ce^{-cz} dz + y_2
= y_2 + b - b e^{-cy_2},
\end{align*}
so $e^{-cy_2} = (y_2 -x_2 +b)/b$ and, therefore, for $z\leq
y_2$,
\begin{align*}
\vphi(z) = ce^{-cz} \geq ce^{-c y_2} = (c/b)(y_2 -x_2 +b).
\end{align*}
We combine this estimate with \eqref{eq:m16.3} and
\eqref{eq:m16.4} to see that, with probability greater than
$p_1$,
\begin{align*}
\bd([x_1,x_2]) &\geq
\frac
{b \int_{y_1}^{y_2} \vphi(z) dz}
{y_2 - y_1 + b \int_{y_1}^{y_2} \vphi(z) dz  }
- (1-d_1)/2 \\
&\geq
\frac
{b \int_{y_1}^{y_2} (c/b)(y_2 -x_2 +b) dz}
{y_2 - y_1 + b \int_{y_1}^{y_2} (c/b)(y_2 -x_2 +b) dz  }
- (1-d_1)/2\\
&=
\frac
{c (y_2-y_1)(y_2 -x_2 +b) }
{y_2 - y_1 + c (y_2-y_1)(y_2 -x_2 +b)  }
- (1-d_1)/2\\
&=
\frac
{c (y_2 -x_2 +b) }
{1 + c (y_2 -x_2 +b)  }
- (1-d_1)/2\\
&\geq
\frac
{c (b -x_2 ) }
{1 + c (b -x_2 )  }
- (1-d_1)/2\\
&\geq 1- (1-d_1)/2
- (1-d_1)/2
= d_1.
\end{align*}
This completes the proof of \eqref{eq:m17.1}. The proof of
\eqref{eq:m17.2} is completely analogous.
\end{proof}

\section{Brownian gas under pressure}
\label{sec:pressure}

In this model, ``particles'' $I^k$ are confined to the interval
$[0,1]$. More precisely, their centers are confined to this
interval. The $k$-th leftmost ``particle'' is represented by an
interval $I^k_t = (X^k_t-\eps/2, X^k_t+\eps/2)$. The intervals
$I^k$ and $I^j$ are always disjoint. The processes $X^k$ are
driven by independent Brownian motions with the diffusion
coefficient $\sigma^2$. When two intervals $I^k$ and $I^j$
collide, they reflect instantaneously. The particles are added
to the system at the left endpoint of $[0,1]$ at a constant
rate. In other words, they are pushed in at the speed $a$, so
that a new particle enters the interval every $\ep/a$ units of
time. As soon as $X^k$ reaches 0, it starts moving as a
Brownian motion reflected at 0. The $k$-th interval is removed
from the system when $X^k$ hits the right endpoint of $[0,1]$.

Formally, we define $\{X^1,X^2,\ldots\}$ to be a collection of
right-continuous, $[0,\infty]$-valued processes such that
  \begin{align}
  &X^k_0=-k\ep+\ep/2 \text{ for all $k$},\\
  &\text{If $S_k=\inf\{t>0: X^k_{t-}=1-\ep/2\}$},\notag\\
  &\qquad\text{then $X^k_t$ is continuous on $[0,S_k)$ and $X^k_t=\infty$ for all $t>S_k$},\\
  &\text{$X^k_t-X^{k+1}_t\ge\ep$ for all $k\ge1$ and all $t\ge0$}, \text{ and}\\
  &dX^k_t = \begin{cases}
    a\,dt &\text{if $t\in[0,k\ep/a)$},\\
    \si\,dB^k_t + dL^k_t - dM^k_t &\text{if $t\in[k\ep/a,S_k)$},
  \end{cases}
  \end{align}
where $a$ and $\si$ are positive constants, $\{B^1,B^2,\ldots\}$ are
iid Brownian motions, and $L^k$ and $M^k$ are nondecreasing processes
such that
  \[
  \int_{k\ep/a}^{S_k} 1_{\{X^k_t-X^{k+1}_t>\ep\}}\,dL^k_t = 0
    \qquad\text{and}\qquad
    \int_{k\ep/a}^{S_k} 1_{\{X^{k-1}_t-X^k_t> \ep\}}\,dM^k_t = 0.
  \]
(Here, we may interpret $X^0\equiv\infty$.)

To construct the solution to this Skorohod problem, consider
first the processes $Y^k_t = X^k_t + k \ep - \ep/2 - at$. These
processes satisfy
  \begin{align}
  &Y^k_0 = 0 \text{ for all $k$},\label{sko1}\\
  &\text{If $S_k=\inf\{t>0: Y^k_{t-}=1-at+(k -1)\ep\}$},\notag\\
  &\qquad\text{then $Y^k_t$ is continuous on $[0,S_k)$ and $Y^k_t=\infty$ for all $t>S_k$},\\
  &\text{$Y^k_t-Y^{k+1}_t\ge0$ for all $k\ge1$ and all $t\ge0$}, \text{ and}\\
  &dY^k_t = \begin{cases}
    0 &\text{if $t\in[0,k\ep/a)$},\\
    \si\,dB^k_t - a\,dt + dL^k_t - dM^k_t &\text{if $t\in[k\ep/a,S_k)$}, \label{sko4}
  \end{cases}
  \end{align}
where
  \[
  \int_{k\ep/a}^{S_k} 1_{\{Y^k_t-Y^{k+1}_t>0\}}\,dL^k_t = 0
    \qquad\text{and}\qquad
    \int_{k\ep/a}^{S_k} 1_{\{Y^{k-1}_t-Y^k_t>0\}}\,dM^k_t = 0.
  \]
Again, we shall construct the processes $\{Y^1,Y^2,\ldots\}$
using order statistics.

Let $Z^k$ be a $[0,\infty)$-valued process, satisfying the SDE
$dZ^k_t = \sigma dB^k_t - a dt$, and reflected at $0$. The
process $Z^k_t$ is defined on the time interval $t\in[k\eps/a,
\infty)$, and starts at $Z^k_{k\eps/a} = 0$. At any time $t \in
[k\eps/a, (k+1)\eps/ a)$, only processes $Z^j$, $1 \leq j \leq
k$, are defined. Let $\lfloor a\rfloor$ denote the greatest
integer less than or equal to $a$, and
\begin{align*}
S_0&=0,\\
\A^1_t &= \{j \in \Z : 1\leq j \leq\lfloor ta/\eps \rfloor\}, \quad t\geq 0,\\
S_1 &= \inf\{t>0: \sup_{j \in \A^1_t} Z^j_t \geq 1- at\},\\
\A^k_t &= \A^{k-1}_t \setminus \{m\in \Z: Z^m_{S_{k-1}} = 1 - at + (k-2) \eps\},
\quad t\geq S_{k-1}, k\geq 2,\\
S_k &= \inf\{t>S_{k-1}:
\sup_{j\in\A^k_t} Z^j_t \geq 1- at+(k-1)\eps\}, \quad k\geq 2.
\end{align*}
Note that it is possible that $\A^k_t = \emptyset$ for some
random $k$ and $t>0$.

\noindent{\it Convention (C)}. For the sake of future reference, it
is convenient to say that the process $Z^m$ is killed at the time
$S_{k-1}$, where $\{m\} = \A^{k-1}_{S_{k-1}} \setminus
\A^{k}_{S_{k-1}}$. In other words, $Z^m$ is killed when the
corresponding interval $I^k$, defined below, hits the right endpoint
of the interval $[0,1]$.

For every $t\in[S_{k-1}, S_k)$, note that there are $\lfloor
ta/\eps \rfloor - (k-1)$ elements in $\A^k_t$. Let $Y^k_t,
Y^{k+1}_t, \dots, Y^{\lfloor ta/\eps \rfloor}_t$ be
reverse-ordered $Z^j_t$'s, $j\in\A^k_t$, that is, $\{Y^k_t,
\dots, Y^{\lfloor ta/\eps \rfloor}_t\} =\{Z^j_t, j\in\A^k_t\}$
and $Y^k_t \ge Y^{k+1}_t\ge \dots \ge Y^{\lfloor ta/\eps
\rfloor}_t$. Let $Y^j_t=\infty$ for $j<k$ and $Y_t^j=0$ for
$j>\lfloor ta/\eps \rfloor$.

It is elementary to check that $\{Y^1,Y^2,\ldots\}$ satisfy
\eqref{sko1}-\eqref{sko4}. We may therefore define
\begin{align*}
X^k_t &= Y^k_t - k\eps
+ \eps/2 + at, \\
I^k_t &= (X^k_t-\eps/2, X^k_t+\eps/2).
\end{align*}
We have to modify slightly the definition \eqref{eq:m17.4} of
density to match the current model. For $t\in[S_{k-1}, S_k)$,
let
\begin{align}\label{eq:m20.1}
\bd([x_1, x_2])= \bd_t([x_1, x_2])=\frac { \lambda\left( [x_1,
x_2] \cap \bigcup_{k \leq j \leq \lfloor ta/\eps \rfloor} I^j_t  \right)}
{x_2 - x_1}.
\end{align}

\begin{theorem}\label{th:m17.5}
Fix arbitrary $0<x_1<x_2<1$, $p_1<1$ and $a,\sigma,c_0>0$.
There exist $t_0< \infty$ and $\eps_0>0$ such that for $t\geq
t_0$ and $\eps \in (0,\eps_0)$,
\begin{align}\label{eq:m17.7}
P\left(\frac {1-x_2} {1-x_2 + \sigma^2/(2a)} -c_0
\leq \bd_t([x_1,x_2])
\leq \frac {1-x_1} {1-x_1 + \sigma^2/(2a)} +c_0
\right) &\geq p_1.
\end{align}
\end{theorem}

Intuitively speaking, the theorem says that the mass density at
$x\in(0,1)$ is close to $(1-x)/(1-x + \sigma^2/(2a))$, for
large $t$ and small $\eps$.

\begin{proof}[Proof of Theorem \ref{th:m17.5}]

We will use the coupling technique. Recall processes\break $Z^1, Z^2,
\dots$ used in the definition of $Y^k$'s---we will use the same
$Z^k$'s to construct auxiliary processes. Fix some $v_1>0$, let
$\wh S_k = \inf\{t\geq 0: Z^k_t = v_1\}$, and let $\wh Z^k_t$
be the process $Z^k$ killed at the time $\wh S_k$. Let $n_t$ be
the number of processes $\wh Z^k$ alive at time $t$. Let $\wh
Y^1_t, \wh Y^2_t, \dots, \wh Y^{n_t}_t$ be ordered $\wh
Z^j_t$'s, that is, $\{\wh Y^1_t, \dots, \wh Y^{n_t}_t\} = \{\wh
Z^j_t, \wh S_j > t\}$ and $\wh Y^1_t \leq \wh Y^2_t\leq \dots
\leq \wh Y^{n_t}_t$. For $t\in[\wh S_{k-1}, \wh S_k)$ and
$j=k,\dots, k + n_t-1$, we let
\begin{align*}
\wh X^j_t &= \wh Y^{n_t + k - j}_t +(n_t + k - j-1) \eps
+ \eps/2 + (t- \lfloor ta/\eps \rfloor \eps /a) a, \\
\wh I^j_t &= (\wh X^j_t-\eps/2, \wh X^j_t+\eps/2).
\end{align*}
Every process $\wh X^j_t$ is defined on the interval
$[j\eps/a,\wh S_{j})$ and it is continuous on this interval.
Although it may not be apparent from the above formulas, the
processes $\wh Y^j, \wh X^j$ and $\wh I^j$ are constructed from
$\wh Z^j$'s in the same way as $ Y^j,  X^j$ and $ I^j$ were
constructed from $ Z^j$'s. We leave the verification of this
claim to the reader.

We will find the Green function $G_{v_1}(v)$ of $\wh Z^k$,
i.e., the density of its occupation measure. Consider a process
$V$ with values in $[-v_1, v_1]$, satisfying the SDE $d V_t =
dB_t - a \sign ( V_t) dt$, where $B$ is Brownian motion, $ V_0
= 0$, and such that $ V$ is killed when it hits $-v_1$ or
$v_1$. Note that the Green function $ G^V_{v_1}(v)$ of $ V$ is
one half of $G_{v_1}(v)$ for $v>0$. The scale function $S(v)$
and the speed measure $m(v)$ for $ V$ can be calculated as
follows (see \citealp{KT2}, pp. 194-195),
\begin{align*}
s(v) & = \exp\left(\int_0^v -(-2 a \sign(x)/\sigma^2) dx\right)
= \exp(2 a v \sign(v)/\sigma^2),\\
S(v) &= \int_0^v s(x) dx
= \frac{\sign(v) \sigma^2}{2a} \left( \exp ( 2 a v \sign(v)/\sigma^2) -1\right),\\
m(v) &= 1/(\sigma^2 s(v)) = (1/\sigma^2) \exp(-2 a v \sign(v)/\sigma^2).
\end{align*}
We will use formula (3.11) on page 197 of \cite{KT2}. In that
formula, we take $x=0$, so $u(0) = 1/2$, by symmetry. We apply
the formula to functions $g(v)$ of the form $g(v) =
\bone_{[v_3, v_4]}(v)$, to conclude that for $v\in(0,v_1)$, the
Green function $ G^V_{v_1}(v)$ is given by
\begin{align*}
 G^V_{v_1}(v)& = (S(v_1) - S(v)) m(v)\\
&= \frac 1 {2a}
(\exp ( 2 a v_1 /\sigma^2) - \exp ( 2 a v /\sigma^2))
\exp(-2 a v /\sigma^2) \\
&= \frac 1 {2a} (\exp ( 2 a (v_1-v) /\sigma^2) -1).
\end{align*}
It follows that
\begin{align*}
G_{v_1}(v) = 2  G^V_{v_1}(v)
= \frac 1 {a} (\exp ( 2 a (v_1-v) /\sigma^2) -1).
\end{align*}

Define $v_0\in(0,\infty)$ by setting
\begin{align*}
\vphi(v) = a G_{v_0}(v)
= \frac 1 {a} (\exp ( 2 a (v_0-v) /\sigma^2) -1),
\end{align*}
and the following condition,
\begin{align}\label{eq:m21.4}
1 &=  \int_0^{v_0} \vphi(v) dv + v_0
= \int _0^{v_0} (\exp ( 2 a (v_0-v) /\sigma^2) -1)dv + v_0\\
&= (\sigma^2/2a) (\exp ( 2 a v_0 /\sigma^2)-1).\nonumber
\end{align}
We define $y_k\in(0,\infty)$, $k=1,2$, by the following
formula,
\begin{align}\nonumber
x_k &= \int_0^{y_k} \vphi(v) dv + y_k\\
&= \int_0^{y_k} (\exp ( 2 a (v_0-v) /\sigma^2) -1) dv + y_k
\nonumber \\
& = (-(\sigma^2/2a)\exp ( 2 a (v_0-v) /\sigma^2) -v)\Big|_{v=0}^{v=y_k}
+ y_k \nonumber \\
&= (\sigma^2/2a)(\exp ( 2 a v_0  /\sigma^2) -
\exp( 2 a (v_0-y_k) /\sigma^2) ). \label{eq:m21.2}
\end{align}

Choose $ y_1 < y_3 < y_4 < y_2$ and $v_1< v_0$ such that,
\begin{align}\label{eq:m21.1}
\frac
{a\int_{y_3}^{y_4} G_{v_1}(v) dv}
{y_2 - y_1 + \int_{y_1}^{y_2} \vphi(z) dz  }
\geq
\frac
{ \int_{y_1}^{y_2} \vphi(z) dz}
{y_2 - y_1 +  \int_{y_1}^{y_2} \vphi(z) dz  }
- c_0.
\end{align}

Recall that $\lceil a\rceil$ denotes the smallest integer
greater than or equal to $a$. Let $\lfloor a\rfloor$ denote the
largest integer smaller than or equal to $a$.

Let $c_1 = 1-p_1$ and $p_2 = 1- c_1/8$. By the law of large
numbers, we can find a large $t_0$ and make $\eps_0>0$ smaller,
if necessary, such that if $t\geq t_0$ and $\eps\in(0,\eps_0)$
then with probability greater than $p_2$, the number of
processes $\wh Z^k_{t}$ in the interval $[0,y_1]$ is smaller
than or equal to $(a/\eps)\int_0^{y_3} G_{v_1}(v) dv$. If this
event holds then there are exactly $\left\lceil
(a/\eps)\int_0^{y_3} G_{v_1}(v) dv\right\rceil$ processes $\wh
Z^k_0$ in some (random) interval $[0,y_5]$ with $y_5 \geq y_1$.
This implies that there are exactly $\left\lceil
(a/\eps)\int_0^{y_3} G_{v_1}(v) dv\right\rceil$ processes $\wh
X^k_0$ in $[0, \eps \left\lceil (a/\eps)\int_0^{y_3} G_{v_1}(v)
dv\right\rceil + y_5]$. For fixed $y_1$ and $y_3$, we make
$v_1< v_0$ larger, if necessary, so that
\begin{align*}
\eps \left\lceil (a/\eps)\int_0^{y_3}
G_{v_1}(v) dv \right\rceil + y_5
&\geq a \int_{0}^{y_3} G_{v_1}(v) dv + y_5
\geq a \int_{0}^{y_3} G_{v_1}(v) dv + y_1 \\
&\geq a\int_{0}^{y_1} G_{v_0}(v) dv + y_1 =
\int_{0}^{y_1} \vphi(v) dv + y_1 =x_1.
\end{align*}
Hence, the number of $\wh X^k_0$'s in the interval $[0, x_1]$
is smaller than or equal to $(a/\eps)\int_0^{y_3} G_{v_1}(v)
dv$, with probability greater than $p_2$.

We can make $t_0$ larger and $\eps_0>0$ smaller, if necessary,
so that by the law of large numbers, if $t\geq t_0$ and
$\eps\in(0,\eps_0)$ then with probability greater than $p_2$,
the number of processes $\wh Z^k_{t_3}$ in the interval
$[0,y_2]$ is greater than or equal to $(a/\eps)\int_0^{y_4}
G_{v_1}(v) dv$. If this event holds then there are exactly
$\left\lfloor (a/\eps)\int_0^{y_4} G_{v_1}(v) dv\right\rfloor$
processes $\wh Z^k_0$ in some (random) interval $[0,y_6]$ with
$y_6 \leq y_2$. This implies that there are exactly
$\left\lfloor (a/\eps)\int_0^{y_4} G_{v_1}(v) dv\right\rfloor$
processes $\wh X^k_0$ in $[0, \eps \left\lceil
(a/\eps)\int_0^{y_4} G_{v_1}(v) dv\right\rceil + y_6]$. Note
that,
\begin{align*}
\eps \left\lfloor (a/\eps)\int_0^{y_4}
G_{v_1}(v) dv \right\rfloor + y_6
&\leq a \int_{0}^{y_4} G_{v_1}(v) dv + y_2
\leq a \int_{0}^{y_2} G_{v_0}(v) dv + y_2
=x_2.
\end{align*}
Hence, the number of $\wh X^k_0$'s in the interval $[0, x_2]$
is greater than or equal to\break $(a/\eps)\int_0^{y_4} G_{v_1}(v)
dv$, with probability greater than $p_2$.

Let $\wh \bd$ be defined as in \eqref{eq:m20.1} but relative to
$\wh I^k$ in place of $I^k$. The two events described in the
last two paragraphs hold simultaneously with probability
greater than $1-c_1/4$. Then the number of $X^k_0$'s in $[x_1,
x_2]$ is greater than or equal to $(a/\eps)\int_{y_3}^{y_4}
G_{v_1}(v) dv$. This and \eqref{eq:m21.1} imply that
\begin{align}\label{eq:m19.3}
\wh\bd_t([x_1,x_2]) &\geq \frac
{\eps\left((a/\eps)\int_{y_5}^{y_6} G_{v_1}(v) dv\right)}
{x_2 - x_1 }
= \frac
{a\int_{y_5}^{y_6} G_{v_1}(v) dv}
{ \int_0^{y_2} \vphi(z) dz + y_2 -  \int_0^{y_1} \vphi(z) dz - y_1 }
\nonumber \\
&= \frac
{a\int_{y_5}^{y_6} G_{v_1}(v) dv}
{y_2 - y_1 + \int_{y_1}^{y_2} \vphi(z) dz  }
\geq
\frac
{ \int_{y_1}^{y_2} \vphi(z) dz}
{y_2 - y_1 +  \int_{y_1}^{y_2} \vphi(z) dz  }
- c_0.
\end{align}

It follows from \eqref{eq:m21.2} that
\begin{align*}
\exp ( 2 a (v_0 -y_2) /\sigma^2) = \exp( 2 a v_0 /\sigma^2)
- 2ax_2 /\sigma^2
\end{align*}
and, therefore, for $v\leq y_2$,
\begin{align*}
\vphi(v) = \exp ( 2 a (v_0-v) /\sigma^2) -1 \geq
\exp ( 2 a (v_0-y_2) /\sigma^2) -1\\
= \exp( 2 a v_0 /\sigma^2)-2ax_2 /\sigma^2 -1.
\end{align*}
We combine this estimate with \eqref{eq:m19.3} to see that,
with probability greater than $1-c_1/4$,
\begin{align}\nonumber
\wh \bd_t([x_1,x_2]) &\geq
\frac
{ \int_{y_1}^{y_2} \vphi(v) dv}
{y_2 - y_1 +  \int_{y_1}^{y_2} \vphi(v) dv  }
- c_0\\
&\geq
\frac
{ \int_{y_1}^{y_2} (\exp( 2 a v_0 /\sigma^2)-2ax_2 /\sigma^2 -1)dv}
{y_2 - y_1 + \int_{y_1}^{y_2} (\exp( 2 a v_0 /\sigma^2)-2ax_2 /\sigma^2 -1)dv  }
- c_0 \nonumber \\
&=
\frac
{ (y_2-y_1)(\exp( 2 a v_0 /\sigma^2)-2ax_2 /\sigma^2 -1) }
{y_2 - y_1 +  (y_2-y_1)(\exp( 2 a v_0 /\sigma^2)-2ax_2 /\sigma^2 -1)  }
- c_0\nonumber \\
&=
\frac
{ \exp( 2 a v_0 /\sigma^2)-2ax_2 /\sigma^2 -1 }
{ \exp( 2 a v_0 /\sigma^2)-2ax_2 /\sigma^2  }
- c_0\nonumber \\
&=
\frac
{ (\sigma^2/2a)(\exp( 2 a v_0 /\sigma^2)-1) -x_2  }
{ (\sigma^2/2a)(\exp( 2 a v_0 /\sigma^2)-1) -x_2 + \sigma^2/2a}
- c_0\nonumber\\
&=
\frac
{ 1 -x_2  }
{ 1 -x_2 + \sigma^2/2a}
- c_0. \label{eq:m21.6}
\end{align}
The last equality follows from \eqref{eq:m21.4}.

Recall that $n_t$ is the number of processes $\wh Z^k$ alive at
time $t$. Note that for any $0\leq t_1< t_2 < \infty$ with $t_2
- t_1 \geq \eps/a$, we have,
\begin{align}\label{eq:m19.1}
n_{t_2} - n_{t_1} \leq (a/\eps) (t_2-t_1).
\end{align}
Fix arbitrary $t_1\geq t_0$ and choose $\delta>0$ such that
\begin{align}\label{eq:m21.3}
 \int_0^{v_1} G_{v_1}(v) dv+\delta
\leq  \int_0^{v_0} G_{v_0}(v) dv.
\end{align}
We make $\eps_0>0$ smaller, if necessary, so that, by the law
of large numbers, if $\eps\in(0,\eps_0)$ then with probability
greater than $1-c_1/4$, for all $s_k$ of the form $s_k = k
\delta/2$, $k=0, \dots, \lfloor t_1/\delta\rfloor +1$, we have
$n_{s_k} < (a/\eps)\left(\int_0^{v_1} G_{v_1}(v)
dv+\delta/2\right)$. It follows from \eqref{eq:m19.1} and
\eqref{eq:m21.3} that for $\eps< a\delta/2$,
\begin{align*}
\sup_{0\leq t \leq t_1}n_{t} <(a/\eps)\left( \int_0^{v_1}
G_{v_1}(v) dv+\delta\right)
\leq (a/\eps) \int_0^{v_0} G_{v_0}(v) dv.
\end{align*}
Suppose that this event holds. Then, for every $t\leq t_1$, the right edge of the rightmost
interval $\wh I^k_t$ is to the left of
\begin{align}\label{eq:m21.5}
\eps(a/\eps) \int_0^{v_0} G_{v_0}(v) dv + v_1
= \int_0^{v_0} G_{v_0}(v) dv + v_0
+ (v_1 -v_0)
= 1 + (v_1- v_0) < 1.
\end{align}
The second equality in the above formula follows from
\eqref{eq:m21.4}.

Recall the definitions given before the statement of the
theorem. A process $Z^k$ is killed when the right end of the
rightmost interval $I^k$ hits 1. Since the processes $\wh Z^k$
are driven by the same Brownian motions as $Z^k$,
\eqref{eq:m21.5} implies that every $Z^k$ has a longer lifetime
than $\wh Z^k$. This implies that $\bd_{t_1}([x_1,x_2]) \geq
\wh \bd_{t_1}([x_1,x_2])$. We combine this with
\eqref{eq:m21.6} to conclude that, with probability greater
than $1-c_1/2$,
\begin{align*}
\bd_{t_1}([x_1,x_2]) \geq \wh \bd_{t_1}([x_1,x_2])
\geq
\frac
{ 1 -x_2  }
{ 1 -x_2 + \sigma^2/2a}
- c_0.
\end{align*}
A completely analogous argument shows that, with probability
greater than $1-c_1/2$,
\begin{align*}
\bd_{t_1}([x_1,x_2]) \leq
\frac
{ 1 -x_1  }
{ 1 -x_1 + \sigma^2/2a}
+ c_0.
\end{align*}
This completes the proof of the theorem.
\end{proof}

\section{Discussion}
\label{sec:disc}

\begin{remark}
In the following remarks we will refer to the model
analyzed in Section~\ref{sec:pressure} as model (C). We will
present another model, which we will call (R). Here, C
represents the ``constant'' rate of influx of new particles,
and R stands for the ``random'' rate of influx. Model (R)
consists of a constant number $n$ of ``particles'' $I^k$ which
are confined to the interval $[0,1]$. The $k$-th leftmost
``particle'' is represented by an interval $I^k_t =
(X^k_t-\eps/2, X^k_t+\eps/2)$. The intervals $I^k$ and $I^j$
are always disjoint. The processes $X^k$ are driven by
independent Brownian motions with the diffusion coefficient
$\sigma^2$. When two intervals $I^k$ and $I^j$ collide, they
reflect instantaneously. The number of particles $n$ is such
that $n\eps=b$, a constant. When $X^k$ hits 1, it jumps to 0.
We conjecture that as $\eps \to 0$, the mass density $\bd$ in
the stationary regime for this process has the density
$(1-x)/(1-x + \sigma^2/(2a))$, just like in model (C), where
$a,\sigma$ and $b$ are related by the following formula,
\begin{align*}
\int_0^1
\frac {1-x} {1-x + \sigma^2/(2a)} = b.
\end{align*}
Heuristically, we expect processes $X^k$ in model (R) to jump
at a more or less constant rate in the stationary regime, so
this is why we believe that models (R) and (C) have the same
hydrodynamic limit. We chose not to analyze model (R) in this
paper as it appears to be harder from the technical point of
view while it seems to illustrate the same phenomenon as model
(C).

\end{remark}

\begin{remark}
Model (R) is closely related to a model studied by T.~Bodineau,
B.~Derrida and J.~Lebowitz (\citealp{BDL}). In their model, one
considers a periodic system of $L$ sites with $N$ particles.
The particles perform random walks but cannot cross each
other---it is the symmetric simple exclusion process. At some
fixed edge, the jump rates are no longer symmetric but jumps
occur with rates $p$ in one direction and $1-p$ in the other
direction. The case $p=1$ corresponds to model (R) described in
the previous remark. In the stationary state, the rescaled
density varies linearly on the unit line segment, with a
discontinuity located where the jump rates are biased. Hence,
away from the singularity, the stationary empirical
distribution is harmonic for the generator of the single
particle process, i.e., Laplacian. This does not apply to the
density of mass $\bd$ in our models (C) and (R).
\end{remark}

\begin{remark}
Since the density $\bd$ of the intervals $I^k$ has the form
$(1-x)/(1-x + \sigma^2/(2a))$, it is elementary to check that
the typical gap size between $I^k$'s is $\eps
\sigma^2/(2a(1-x))$. In a model with infinitely small particles
$X^k$, the gap size is also $c/(1-x)$ but we do not have any
heuristic explanation why the two functions representing the
typical gap size should have the same form in both models.
\end{remark}

\begin{remark}\label{rem:m21.1}
We conjecture that the motion of an individual tagged particle
$I^k$ in model (C) converges, as $\eps\to 0$, to a
deterministic motion with fluctuations having the ``fractional
Brownian motion'' structure. In other words, we conjecture that
the fluctuations are Gaussian with the local scaling of space
and time given by $\Delta x = (\Delta t)^{1/4}$. Our conjecture
is inspired by the results in \cite{DGL,H,S,Sw08} on families
of one dimensional Brownian motions reflecting from each other.
\end{remark}

\begin{remark}
A $d$-dimensional counterpart of model (C) can be represented
as follows. Let $I^k$ be balls with radius $\eps$ and center
$X^k$. Our $d$-dimensional model consists of a constant number
$n$ of $I^k$'s which are confined to the cube $[0,1]^d$. The
balls $I^k$ and $I^j$ are always disjoint. The processes $X^k$
are driven by independent $d$-dimensional Brownian motions with
the diffusion coefficient $\sigma^2$. When two balls $I^k$ and
$I^j$ collide, they reflect instantaneously. Let $S_\ell$ and
$S_r$ be two opposite $(d-1)$-dimensional sides on the boundary
of $[0,1]^d$. Balls $I^k$ are pushed into the cube through
$S_\ell$ at a constant rate $a$, i.e., $\eps^{-(d-1)}$ balls
are pushed into the cube every $\eps/a$ units of time,
uniformly over $S_\ell$. Once inside the cube, the balls
reflect from all sides except $S_r$. When a ball hits $S_r$, it
is removed from the cube. We conjecture that in the stationary
regime, when $\eps$ is small, the density of the mass analogous
to $\bd$ will be a function of the distance $x$ from $S_\ell$,
i.e., a function of depending only on one coordinate. We do not
see any obvious reason why the density should have the form
$(1-x)/(1-x + \sigma^2/(2a))$. In relation to Remark
\ref{rem:m21.1}, we conjecture that the motion of a tagged
particle in the present model is diffusive, with the diffusion
coefficient depending on $x$. If this is true, it means that
the ``pressure'' applied to particles in one direction can have
a dampening effect on the size of oscillations of an individual
particle in orthogonal directions.

\end{remark}

\bibliographystyle{alea2}
\bibliography{07-10}

\end{document}